\documentclass[lettersize,journal]{IEEEtran}
\usepackage[justification=centering]{caption}
\usepackage{booktabs}
\usepackage{amsmath,amsfonts}
\usepackage{cases}
\usepackage{array}
\usepackage{textcomp}

\usepackage{multirow}
\usepackage{pbox}
\usepackage{makecell}

\usepackage{stfloats}
\usepackage{url}
\usepackage{verbatim}
\usepackage{epstopdf}
\usepackage{graphicx}
\usepackage{subfigure}
\usepackage[justification=centering]{caption}
\usepackage{amsmath}
\usepackage{color}
\usepackage{algorithm}
\usepackage{algpseudocode}

\usepackage{cite}

\newtheorem{remark}{Remark}
\newtheorem{theorem}{Theorem}
\newtheorem{lemma}{Lemma}

\newtheorem{proof*}{\noindent Proof}

\begin{document}
\title{Optimal control of discrete-time nonlinear systems}
\author{
	\vskip 1em
	Chuanzhi Lv, Xunmin Yin, Hongdan Li, Huanshui Zhang, \emph{Senior Member, IEEE}

	\thanks{This work was supported by the 
		Original Exploratory Program Project of National Natural Science Foundation of China (62450004),  
		the Joint Funds of the National Natural Science Foundation of China (U23A20325),  
		the Major Basic Research of Natural Science Foundation of Shandong Province (ZR2021ZD14), 
		the National Natural Science Foundation of China (62103241),
		the Shandong Provincial Natural Science Foundation (ZR2021QF107), 
		and the Development Plan for Youth Innovation Teams in Higher Education Institutions in Shandong Province (2022KJ222). (\textit{Corresponding author: Huanshui Zhang.})
		
		C. Lv, X. Yin and H. Li are with the College of Electrical Engineering and Automation, Shandong University of Science and Technology, Qingdao 266590, China (e-mail: lczbyq@sdust.edu.cn; yxm\_90711@sdust.edu.cn; lhd2022@sdust.edu.cn).
		
		H. Zhang is with the College of Electrical Engineering and Automation, Shandong University of Science and Technology, Qingdao 266590, China, and  also with the School of Control Science and Engineering, Shandong University, Jinan, Shandong 250061, China (e-mail: hszhang@sdu.edu.cn).
		
		This work has been submitted to the IEEE for possible publication. Copyright may be transferred without notice, after which this version may no longer be accessible.
	}
}

\maketitle
	
\begin{abstract}
	This paper focuses on optimal control problem for a class of discrete-time nonlinear systems. 
	In practical applications, computation time is a crucial consideration when solving nonlinear optimal control problems, especially under real-time constraints. While linearization methods are computationally efficient, their inherent low accuracy can compromise control precision and overall performance.
%
	To address this challenge, this study proposes a novel approach based on the optimal control method. 
	Firstly, the original optimal control problem is transformed into an equivalent optimization problem, which is resolved using the Pontryagin's maximum principle, and a superlinear convergence algorithm is presented. 
	Furthermore, to improve computation efficiency, explicit formulas for computing both the gradient and hessian matrix of the cost function are proposed. 
	Finally, the effectiveness of the proposed algorithm is validated through simulations and experiments on a linear quadratic regulator problem and an automatic guided vehicle trajectory tracking problem, demonstrating its ability for real-time online precise control.
\end{abstract}


\markboth{}%
{}

\definecolor{limegreen}{rgb}{0.2, 0.8, 0.2}
\definecolor{forestgreen}{rgb}{0.13, 0.55, 0.13}
\definecolor{greenhtml}{rgb}{0.0, 0.5, 0.0}

\section{Introduction}
\IEEEPARstart{O}{ptimal} control problem (OCP) is fundamentally an optimization problem (OP) focused on determining the control inputs for a dynamical system to minimize or maximize a cost function over time. It can be argued that many natural processes are essentially driven by optimization, reflecting the universal principle of optimality \cite{liberzon2011calculus}. In a world increasingly shaped by optimization-based control, the growing diversity of applications demonstrates its essential role and wide-ranging impact across various fields. For example, the optimization-based control strategies have enabled bipedal, quadrupedal, and humanoid robots to walk reliably in standard environments, leading to their first practical deployments, such as Boston Dynamics' Spot robot \cite{wensing2023optimization}. 
It is important to note that most real-world systems exhibit inherent nonlinearities. While linearization is commonly employed in control system design for simplification, it often results in the loss of important nonlinear characteristics \cite{cimen2004nonlinear, ren2020optimal}. 
As the demand for enhanced performance increases, it becomes increasingly critical to account for the nonlinear nature of the system in the design and analysis process.

\par Viewed from a static perspective, the nonlinear OCP can be solved using a variety of existing optimization algorithms \cite{boyd2004convex}. Model predictive control (MPC) represents a practical implementation of this theory, effectively transforming the complex OCP into a series of smaller optimization problems that are solved in real-time \cite{kang2022tracking}. In autonomous driving scenarios, MPC was used to optimize the angular velocity and steering wheel angles of automatic guided vehicles (AGV) for improving tracking accuracy \cite{deng2023model}. In \cite{zhao2024nonlinear}, the control strategy for AGV parallel parking was built around MPC, which ensured smooth transitions, stability, and effective handling of non-convex state constraints. To reduce the computational demands of real-time optimization, an event-triggered MPC approach was applied to study path tracking control in AGV \cite{zhou2023event}. Although there are many available optimization algorithms, the complexity of the optimization problem often increases with longer terminal time, leading to a greater computational burden.

\par Conversely, considering a dynamic perspective, the Pontryagin's maximum principle is the key mathematical tool for solving nonlinear OCP, focusing on the solution of highly coupled nonlinear forward and backward difference equations (FBDEs), which reveal the essence of optimal control. In \cite{wang2024decentralized}, the decentralized optimal control and stabilization problems for linearly interconnected systems with asymmetric information were effectively addressed by solving stochastic FBDEs. Unfortunately, analytical solutions to nonlinear FBDEs are nearly impossible to obtain, except in the case of linear systems. A well-known numerical method for solving nonlinear FBDEs is the
method of successive approximation (MSA) \cite{chernousko1982method}. In the field of machine learning, \cite{chen2022self} proposed a self-healing neural network architecture that employed this method to automatically detect and fix the vulnerability issue, thereby enhancing the robustness of the neural network. However, the direct application of MSA, a gradient-based method, often results in a slow convergence rate during computation. 

\par To overcome the limitations of the aforementioned methods, \cite{teo2021applied} introduced an approach that combines nonlinear optimization with the maximum principle to improve the efficiency of solving the OCP. This method primarily employs FBDEs for gradient computation, replacing the traditional gradient computation in optimization. In \cite{li20223D}, the hybrid gradient-based method was applied to simultaneously optimize unmanned aerial vehicle trajectories and control signals, minimizing mission completion time and energy consumption while ensuring smooth and practical implementation. Similarly, the control gains of a controller based on a fully actuated system approach were optimized using the hybrid method to ensure high-precision trajectory tracking while satisfying angular and actuator constraints \cite{tian2023high}. Besides, the approach has also been extended to address the OCP of switched nonlinear systems \cite{zhu2022sequential,sun2023numerical}. It should be noted that both the selection of optimization algorithms and gradient computation are crucial for computational efficiency, with each playing an indispensable role. However, the majority of existing studies devote relatively less attention to the design of optimization algorithms. Furthermore, the current methods for gradient calculation are often not readily applicable to higher-order optimization algorithms.

\par In this paper, the OCP for discrete-time nonlinear systems is considered. 
Unlike previous research that relies on existing optimization algorithms \cite{teo2021applied,li20223D,tian2023high,zhu2022sequential,sun2023numerical}, this work presents a novel optimization approach from the perspective of optimal control, featuring a superlinear convergence property.
Following the work of \cite{zhang2024optimization}, this study specifically addresses dynamic problems, with a particular focus on the challenges associated with effectively computing the gradient and hessian matrix.
To tackle these challenges, a new approach based on the maximum principle is introduced for the first time, transforming the solution process into solving two sets of FBDEs.
This provides a unified computational approach for developing second-order optimal control algorithms.
The contributions of this paper can be summarized as follows:
\begin{itemize}
	\setlength{\itemindent}{1.5mm}
	\item[1)]The discrete-time nonlinear OCP is transformed into an OP. By constructing a novel OCP, a rapidly convergent iterative algorithm for solving the OP is presented, establishing a new framework for solving the OCP based on the principle of optimal control theory.
	\item[2)]A unified and explicit formula based on the maximum principle for computing the gradient and hessian matrix in OCP is proposed in this work, which reduces the computational burden. Furthermore, this method can be generalized to other second-order algorithms.
	\item[3)]The proposed optimal control algorithm is validated numerically and experimentally on linear quadratic regulator (LQR) and AGV systems, with results demonstrating its capability to achieve real-time online precise control computation.
\end{itemize}

\par The rest of this paper is organized as follows. In section II, the discrete-time nonlinear OCP is given. 
In section III, the control design and the iterative algorithm's implementation and details are presented. Section IV presents the simulation and experimental validation for both LQR and AGV systems. Section V concludes the work.

\par \textit{Notations:} $\mathbb{R}^n$ denotes the n-dimensional Euclidean space; $A^{T}$ represents the transpose of matrix $A$; $\frac{\partial p}{\partial q}$ denotes partial derivative of function $p$ with respect to $q$. $\textbf{0}$ denotes the zero vector.

\section{Problem Statement}
Consider the following discrete-time nonlinear OCP
\begin{alignat}{1}
	&\underset{u_{[0,N]}}{\text{min}} \  J = \sum_{k=0}^{N}\Phi(x_k,u_k), \label{OCP1_cost}\\
	&\hspace*{0.66em}  \mbox{s.t.}\ \hspace*{0.30em} x_{k+1} = F(x_k, u_k), \label{OCP1_sys}  \\
	\nonumber & \hspace*{2.47em} x_0 = \sigma,\ k = 0,1,...,N,
\end{alignat}
in which $\Phi\in\mathbb{R}$ is the running cost, $F\in\mathbb{R}^n$ is a given deterministic function, $x_k \in \mathbb{R}^n$ is the state variable, $u_k \in \mathbb{R}^m$ is the control variable, and $x_0$ represents the given initial state. It should be noted that the cost function (\ref{OCP1_cost}) considered here is in the lagrange form, which encompasses the bolza form as a special case when $\Phi(x_N,u_N)$ is set to $\psi(x_N)$ at $k=N$ \cite{naidu2002optimal}.

\par Notice that the nonlinear dynamical system (\ref{OCP1_sys}) can be reformulated as
\begin{align}
	\begin{bmatrix}
		x_1\\
		x_2\\
		\vdots\\
		x_N
	\end{bmatrix}
	=
	\begin{bmatrix}
		F(x_0,u_0)\\
		F(F(x_0,u_0),u_1)\\
		\vdots\\
		F(F(x_{N-2},u_{N-2}),u_{N-1})
	\end{bmatrix}, 
\end{align} 
where the state $x_k$ can be regarded as a function of $(x_0,u_0,u_1,$ $...,u_{k-1})$. For convenience, denote $z=[u_0^{T},u_1^{T},...,u_{N}^{T}]^{T}$. Then the OCP (\ref{OCP1_cost})-(\ref{OCP1_sys}) can be restated as an OP
\begin{align}
	&\underset{z}{\text{min}}\ J(x_0,z). \label{OP1_cost}
\end{align}
\begin{remark}
	The OCP (\ref{OCP1_cost})-(\ref{OCP1_sys}) is equivalent to the OP (\ref{OP1_cost}). From the perspective of optimal control, the OCP (\ref{OCP1_cost})-(\ref{OCP1_sys}) can be solved using the Pontryagin's maximum principle, which involves solving highly coupled FBDEs. Additionally, from the perspective of nonlinear optimization, the OP (\ref{OP1_cost}) can be tackled using existing optimization method, such as gradient descent algorithm. In the following, the strengths of the aforementioned methods will be combined to present an efficient iterative numerical algorithm.
\end{remark}

\section{Main Results}
In this section, in order to solve the OP (\ref{OP1_cost}), the problem is reformulated as a new OCP from an optimal control perspective, and a rapidly convergent iteration algorithm is presented. Subsequently, the explicit formulas for computing the gradient and hessian matrix of the cost function are proposed by solving the corresponding FBDEs.
\subsection{Optimal Control Design}
In order to solve the optimal controller $u_k$ in (\ref{OCP1_sys}), rather than directly solving the OCP (\ref{OCP1_cost})-(\ref{OCP1_sys}), the focus is placed on its equivalent OP (\ref{OP1_cost}). In what follows, the OP (\ref{OP1_cost}) is reconsidered from a new perspective, namely that of optimal control.

\par Define the following new OCP:
\begin{alignat}{1}
	&\underset{\upsilon_{[0,M-1]}}{\text{min}} \hspace*{-0.2em} J_M \hspace*{-0.2em}=\hspace*{-0.3em} \sum_{i=0}^{M-1}[J(x_0,z_i) + \frac{1}{2}\upsilon_i^TR\upsilon_i] \hspace*{-0.05em}+\hspace*{-0.05em} J(x_0,z_M), \label{OCP3_cost}\\
	&\hspace*{0.73em}  \mbox{s.t.} \hspace*{0.93em} z_{i+1} = z_i + \upsilon_i, \label{OCP3_sys}  \\
	\nonumber & \hspace*{2.83em} z_0 = \eta,\ i = 0,1,...,M-1, 
\end{alignat}
where $z_i\in\mathbb{R}^{m(N+1)}$ and $\upsilon_i\in\mathbb{R}^{m(N+1)}$ are the state and control variables, respectively, $z_0$ is the initial state value, $R$ is the positive definite matrix. It should be noted that the design of the cost function (\ref{OCP3_cost}) is based on that of OP (\ref{OP1_cost}). 
At each time step, the controller $v_i$ minimizes the sum of costs over future steps, i.e., $\sum_{j=i+1}^{M}J(x_0,z_j)$, thereby accelerating the convergence of the optimization algorithm.
Furthermore, by integrating control energy minimization into the loss function (\ref{OCP3_cost}), the algorithm can achieve enhanced stability.


\par In order to solve the OCP (\ref{OCP3_cost})-(\ref{OCP3_sys}), the maximum principle is applied to derive the following equations:
\begin{align}
	\gamma_i &= \gamma_{i+1} + \nabla J(x_0,z_i),  \label{gamma} \\  
	0 &= R\upsilon_i + \gamma_i, \label{equilibrium}
\end{align}
where (\ref{gamma}) is the costate equation with the terminal condition $\gamma_M=\nabla J(x_0,z_M)$, and (\ref{equilibrium}) is the equilibrium condition. It is known that any solution to the OCP must satisfy the FBDEs. Therefore, the optimal control and state of the \hspace*{-0.05em}OCP\hspace*{-0.05em} (\ref{OCP3_cost})-(\ref{OCP3_sys}) can be obtained by solving the FBDEs (\ref{OCP3_sys})-(\ref{equilibrium}). It is evident from the FBDEs (\ref{OCP3_sys})-(\ref{equilibrium}) that when the state $z_i$ in the system (\ref{OCP3_sys}) stabilizes, the terminal state $z_M$ represents the optimal solution to the OP (\ref{OP1_cost}). 

\par Although FBDEs (\ref{OCP3_sys})-(\ref{equilibrium}) can be numerically solved using the existing MSA, its convergence rate is relatively slow. To improve computational efficiency, an iterative numerical algorithm has been developed by combining the Taylor expansion with mathematical induction in \cite{zhang2024optimization}, as outlined below.
\begin{align}
	z_{i+1} &\hspace*{-0.1em}=\hspace*{-0.1em} z_i - g_i(z_i), \label{z_i+1_2}
	\\
	g_j(z_i) &\hspace*{-0.1em}=\hspace*{-0.1em}(R \hspace*{-0.1em}+\hspace*{-0.1em} \nabla^2 J(x_0,z_i))^{-1}(\nabla J(x_0,z_i) \hspace*{-0.2em}+\hspace*{-0.2em} Rg_{j-1}(z_i)), \label{g_i_2}
	\\
	g_0(z_i) &\hspace*{-0.1em}=\hspace*{-0.1em} (R \hspace*{-0.1em}+\hspace*{-0.1em} \nabla^2 J(x_0,z_i))^{-1}\nabla J(x_0,z_i), \label{g_M_1_2}
\end{align}
where $i=0,1,...,M-1$ and $j=1,...,i$. 
\begin{remark}
	In the iterative algorithm (\ref{z_i+1_2})-(\ref{g_M_1_2}), the update step $g_i$ requires inner loop iterations starting from $g_0$. 
	It should be noted that the algorithm faces considerable challenges in handling dynamic problem, i.e., the OCP (\ref{OCP1_cost})-(\ref{OCP1_sys}), due to the computational complexity involved in evaluating $\nabla J$ and $\nabla^2 J$. 
	The following subsections provide the explicit iterative solution algorithms. For details on the convergence analysis and the superlinear convergence rate of the algorithm (\ref{z_i+1_2})-(\ref{g_M_1_2}), please refer to \cite{zhang2024optimization}.
\end{remark}

\subsection{Computation of Gradient}
From the OCP (\ref{OCP1_cost})-(\ref{OCP1_sys}), it can be seen that the cost function $J$ is a complex function on the control variables $u_k$ to be optimized. To enhance computational efficiency, the gradient formula for the cost function $J$ is provided in this subsection.

\par To simplify presentation, a one-dimensional case for the control variables $u_k$ is considered, with a natural extension to higher dimensions. Define the following Hamiltonian 
\begin{align}
	H(x_k,u_k,\lambda_{k+1}) = \Phi(x_k,u_k) + \lambda_{k+1}^{T}F(x_k,u_k), \label{H_N}
\end{align}
where $k=0,1,...,N$, and $\lambda_k\in\mathbb{R}^n$ is the costate variable.
Applying the Pontryagin's maximum principle to the OCP (\ref{OCP1_cost})-(\ref{OCP1_sys}), the main result is summarized in the following lemma.
\begin{lemma}
	For the OP (\ref{OP1_cost}), the gradient of $J$ with respect to $z$ is given by
	\begin{align}
		\nabla J = \left[\frac{\partial H(x_0,u_0,\lambda_{1})}{\partial u_0},...,\frac{\partial H(x_{N},u_{N},\lambda_{N+1})}{\partial u_{N}}\right]^T, \label{gradient}
	\end{align} 
	where $x_k$ and $H$ are given by (\ref{OCP1_sys}) and (\ref{H_N}), and $\lambda_k$ is given by 
	\begin{align}
		\lambda_k = \frac{\partial H(x_k,u_k,\lambda_{k+1})}{\partial x_k}, \label{lambda_k}
	\end{align}
	with the terminal costate value $\lambda_{N+1} = \textbf{0}$.
\end{lemma}

\par \textit{Proof:} For the proof, please refer to \cite{naidu2002optimal}, which is omitted here.

\par The implementation steps of Lemma 1 are summarized as follows:

\par \textit{Step 1:} For the given control $u_k, k=0,1,...,N$, compute the forward difference equation (\ref{OCP1_sys}) from $k=0$ to $k=N-1$ to obtain the state $x_k, k=0,1,...,N$ using the given initial state value $x_0$.

\par \textit{Step 2:} Compute the backward difference equation (\ref{lambda_k}) from $k=N$ to $k=1$ to obtain the costate $\lambda_k, k=1,...,N+1$ using the given terminal costate value $\lambda_{N+1}$.

\par \textit{Step 3:} Compute the gradient of $J$ according to (\ref{gradient}).

\subsection{Computation of Hessian Matrix}
In this subsection, the second derivative of the cost function $J$ is computed based on Lemma 1, and the explicit formula for hessian matrix computation will be presented.

\par During the computation process, we solve for the elements of the hessian matrix line by line. For convenience, denote 
\begin{align}
	\nonumber G(x_k,u_k,\lambda_{k+1})=\frac{\partial H(x_k,u_k,\lambda_{k+1})}{\partial x_k},\\
	\nonumber L(x_k,u_k,\lambda_{k+1})=\frac{\partial H(x_k,u_k,\lambda_{k+1})}{\partial u_k}.
\end{align}
To compute the elements in the $i$th row of the hessian matrix, construct the following OCP:
\begin{alignat}{1}
	&\underset{u_{[0,N]}}{\text{min}} \hspace*{1.2em}  J_i = L(x_i,u_i,\lambda_{i+1}), \label{OCP2_cost}\\
	&\hspace*{0.66em}  \mbox{s.t.}\ \hspace*{1.15em} x_{k+1} = F(x_k, u_k),\ k = 0,1,...,N-1,  \label{OCP2_forward} \\
	& \hspace*{3.34em} \lambda_k = G(x_k,u_k,\lambda_{k+1}),\ k = N,...,1, \label{OCP2_backward}\\
	\nonumber & \hspace*{3.34em} x_0=\alpha,\ \lambda_{N+1} = \textbf{0},
\end{alignat}
where $i=0,1,...,N$, and the cost function (\ref{OCP2_cost}) is subject to constraints imposed by the FBDEs (\ref{OCP2_forward})-(\ref{OCP2_backward}).

\par Let $L(x_i,u_i,\lambda_{i+1})=\sum_{k=0}^{N}L(x_k,u_k,\lambda_{k+1})$, where $L(x_k,u_k,\lambda_{k+1})=0, k\neq i$. Define the Hamiltonian as follows when $k=i$
\begin{align}
	\nonumber &H_i(x_k,u_k,\lambda_{k+1},\alpha_{k+1},\beta_{k})\\
	\nonumber=&L(x_k,u_k,\lambda_{k+1}) + \alpha_{k+1}^TF(x_k,u_k)  \\ &\hspace*{6.55em}+\beta_{k}^TG(x_k,u_k,\lambda_{k+1}). \label{H_i_L}
\end{align}
Otherwise, it is given by
\begin{align}
	\nonumber &H_i(x_k,u_k,\lambda_{k+1},\alpha_{k+1},\beta_{k})\\
	=&\alpha_{k+1}^TF(x_k,u_k) + \beta_{k}^TG(x_k,u_k,\lambda_{k+1}). \label{H_i}
\end{align}
where $\alpha$ and $\beta$ are the Lagrange multiplier, $i=0,1,...,N$, and  $k=0,1,...,N$.
Besides, define the following FBDEs based on the Hamiltonian (\ref{H_i_L})-(\ref{H_i})
\begin{align}
	\beta_{k+1} &= \frac{\partial H_i(x_k,u_k,\lambda_{k+1},\alpha_{k+1},\beta_{k})}{\partial \lambda_{k+1}}, \label{beta}\\
	\alpha_k &= \frac{\partial H_i(x_k,u_k,\lambda_{k+1},\alpha_{k+1},\beta_{k})}{\partial x_k}, \label{alpha}
\end{align}
with $\beta_0=\textbf{0}$, $\alpha_{N+1}=\textbf{0}$.

\par The main result of this subsection is given below.
\begin{theorem}
	For the OP (\ref{OP1_cost}), the elements of the $i$th row of the hessian matrix of the cost function $J$ with respect to $z$ are given by
	\begin{align}
		\nabla^2J(i) = \left[\frac{\partial H_0^{i}}{\partial u_0},\frac{\partial H_1^{i}}{\partial u_1},...,\frac{\partial H_{N}^i}{\partial u_{N}}\right], \label{hessian}
	\end{align}
	where $H^i_k=H_i(x_k,u_k,\lambda_{k+1},\alpha_{k+1},\beta_{k})$, and $x_k, \lambda_{k+1}, \beta_{k}$ and $\alpha_{k+1}$ are given by (\ref{OCP2_forward}), (\ref{OCP2_backward}), (\ref{beta}) and (\ref{alpha}), respectively.
\end{theorem}

\par \textit{Proof:} From the OCP (\ref{OCP2_cost})-(\ref{OCP2_backward}), we formulate an optimal augmented cost function by adjoining the original cost function (\ref{OCP2_cost}) with the FBDEs (\ref{OCP2_forward})-(\ref{OCP2_backward}) using the Lagrange multiplier $\alpha_{k+1}$ and $\beta_{k}$ as 
\begin{align}
	\nonumber J_i^* = \sum_{k=0}^{N}\big\{L(x_k^*,u_k^*,\lambda_{k+1}^*) &+ \alpha_{k+1}^{*T}[F(x_k,u_k)-x_{k+1}^*] \\ &\hspace*{-2.3em}+\beta_{k}^{*T}[G(x_k^*,u_k^*,\lambda_{k+1}^*)-\lambda_k^*]\big\}. \label{Opti_cost}
\end{align}
Let $x_k$, $u_k$, $\lambda_{k+1}$, $x_{k+1}$, and  $\lambda_k$ deviate from their optimal values $x_k^*$, $u_k^*$, $\lambda_{k+1}^*$, $x_{k+1}^*$, and  $\lambda_k^*$ by variations $\delta x_k$, $\delta u_k$, $\delta \lambda_{k+1}$, $\delta x_{k+1}$, and $\delta \lambda_k$, respectively, which means that 
\begin{align}
	\nonumber \xi_k &= \xi_k^* + \delta \xi_k,\\
	\nonumber \zeta_{k+1} &= \zeta_{k+1}^* + \delta \zeta_{k+1}.
\end{align}
where $\xi = x,u,\lambda$ and $\zeta = \lambda, x$. Based on these variations, the cost function (\ref{Opti_cost}) becomes
\begin{align}
	\nonumber J_i = \sum_{k=0}^{N-1}\big\{L(x_k,u_k,\lambda_{k+1}) &+ \alpha_{k+1}^{T}[F(x_k,u_k)-x_{k+1}] \\ &\hspace*{-2.3em}+\beta_{k}^{T}[G(x_k,u_k,\lambda_{k+1})-\lambda_k]\big\}. \label{NotOpti_cost}
\end{align}

\par It is well known that the necessary condition for optimal control is that the first variation must be equal to zero. Therefore, by applying the Taylor series expansion to (\ref{NotOpti_cost}) along with (\ref{Opti_cost}), it follows that 
\begin{align}
	\nonumber \delta J_i=\sum_{k=0}^{N}\big[\delta x_k^T\frac{\partial V^i_k}{\partial x_k^*} + \delta u_k^T\frac{\partial V^i_k}{\partial u_k^*} + \delta \lambda_{k+1}^T\frac{\partial V^i_k}{\partial \lambda_{k+1}^*}  \\
	+\delta x_{k+1}^T\frac{\partial V^i_k}{\partial x_{k+1}^*} +  \delta \lambda_k^T\frac{\partial V^i_k}{\partial \lambda_k^*}\big], \label{first_varia}
\end{align}
where 
\begin{align}
	\nonumber V_k^i &= V_i(x_k^*,u_k^*,\lambda_{k+1}^*,\alpha_{k+1}^*,x_{k+1}^*,\beta_{k}^*,\lambda_k^*)\\
	\nonumber &=L(x_k^*,u_k^*,\lambda_{k+1}^*) + \alpha_{k+1}^{*T}[F(x_k,u_k)-x_{k+1}^*] \\ \nonumber&\hspace*{8.1em}+\beta_{k}^{*T}[G(x_k^*,u_k^*,\lambda_{k+1}^*)-\lambda_k^*]. 
\end{align}
To express the variation $\delta x_{k+1}$ in terms of $\delta x_k$, examine the fourth expression in (\ref{first_varia}) below.
\begin{align}
	\nonumber&\sum_{k=0}^{N}\delta x_{k+1}^T\frac{\partial V^i_k}{\partial x_{k+1}^*} = \sum_{k=1}^{N+1}\delta x_k^T\frac{\partial V^i_{k-1}}{\partial x_k^*}\\
	= &\sum_{k=0}^{N}\delta x_k^T\frac{\partial V^i_{k-1}}{\partial x_k^*} + (\delta x_{N+1}^T\frac{\partial V^i_{N}}{\partial x_{N+1}^*} - \delta x_0^T\frac{\partial V^i_{-1}}{\partial x_0^*}). \label{varia_x}
\end{align}
Similarly, consider the fifth expression in (\ref{first_varia}) below.
\begin{align}
	\nonumber&\sum_{k=0}^{N}\delta \lambda_k^T\frac{\partial V^i_k}{\partial \lambda_k^*} = \sum_{k=-1}^{N-1}\delta \lambda_{k+1}^T\frac{\partial V^i_{k+1}}{\partial \lambda_{k+1}^*}\\
	= &\sum_{k=0}^{N}\delta \lambda_{k+1}^T\frac{\partial V^i_{k+1}}{\partial \lambda_{k+1}^*} + (\delta \lambda_0^T\frac{\partial V^i_0}{\partial \lambda_0^*} - \delta \lambda_{N+1}^T\frac{\partial V^i_{N+1}}{\partial \lambda_{N+1}^*}). \label{varia_lam} 
\end{align}
Substituting the expressions (\ref{varia_x})-(\ref{varia_lam}) into (\ref{first_varia}), the first variation $\delta J_i$ can be reformulate as 
\begin{align}
	\nonumber \delta J_i=\sum_{k=0}^{N}&\big[\delta x_k^T(\frac{\partial V^i_k}{\partial x_k^*} + \frac{\partial V^i_{k-1}}{\partial x_k^*}) + \delta \lambda_{k+1}^T(\frac{\partial V^i_k}{\partial \lambda_{k+1}^*} + \frac{\partial V^i_{k+1}}{\partial \lambda_{k+1}^*})  \\
	\nonumber &+ \delta u_k^T\frac{\partial V^i_k}{\partial u_k^*}\big] + (\delta x_{N+1}^T\frac{\partial V^i_{N}}{\partial x_{N+1}^*} - \delta x_0^T\frac{\partial V^i_{-1}}{\partial x_0^*}) \\
	&+(\delta \lambda_0^T\frac{\partial V^i_0}{\partial \lambda_0^*} - \delta \lambda_{N+1}^T\frac{\partial V^i_{N+1}}{\partial \lambda_{N+1}^*}). \label{first_varia_ed}
\end{align}
Combining the Hamiltonian (\ref{H_i_L}) and (\ref{H_i}), we have 
\begin{align}
	\nonumber V_k^i = H_i(x_k^*,u_k^*,\lambda_{k+1}^*,\alpha_{k+1}^*,\beta_{k}^*) - \alpha_{k+1}^{*T}x^*_{k+1} - \beta_{k}^{*T}\lambda_k. 
\end{align}
Taking into account the arbitrary nature of the variations in (\ref{first_varia_ed}), we can obtain 
\begin{align}
	\frac{\partial H_i(x_k^*,u_k^*,\lambda_{k+1}^*,\alpha_{k+1}^*,\beta_{k}^*)}{\partial u_k^*} = 0,  \label{euqili_condi}
\end{align}
where $x_k^*, \lambda_{k+1}^*, \beta_{k}^*$ and $\alpha_{k+1}^*$ satisfy (\ref{OCP2_forward}), (\ref{OCP2_backward}), (\ref{beta}) and (\ref{alpha}), respectively. When the given $u_k$ does note meet the equilibrium condition (\ref{euqili_condi}), the left side of equation (\ref{euqili_condi}) represents the $k$th element of the gradient, i.e., the (\ref{hessian}) represents the $i$th row of the hessian matrix of the cost function $J$. The proof is complete.

\par Building on the three steps of the gradient computation, the implementation steps of Theorem 1 are summarized as follows:

\par \textit{Step 1:} For the given $x_k$, $u_k$, and $\lambda_k$, compute the forward difference equation (\ref{beta}) from $k=0$ to $k=N-1$ to obtain the $\beta_k, k=0,1,...,N$ using the given initial value $\beta_0$.

\par \textit{Step 2:} Compute the backward difference equation (\ref{alpha}) from $k=N$ to $k=1$ to obtain the $\alpha_k, k=1,...,N+1$ using the given terminal value $\alpha_{N+1}$.

\par \textit{Step 3:} Compute the elements of the $i$th row of the hessian matrix of $J$ according to (\ref{hessian}).

\par \textit{Step 4:} Repeat steps 1 to 3 until all elements of the hessian matrix are computed.

\begin{remark}
	The maximum principle is applied three times throughout this work. In the first application, it solves OCP (\ref{OCP3_cost})-(\ref{OCP3_sys}) to develop the iterative optimization algorithm for OP (\ref{OP1_cost}). In the second, it is used to derive the gradient of the cost function $J$ with respect to $z$ by solving OCP (\ref{OCP1_cost})-(\ref{OCP1_sys}). The final application involves solving OCP (\ref{OCP2_cost})-(\ref{OCP2_backward}) to provide the explicit formula for the second-order derivative of the cost function $J$ with respect to $z$.
\end{remark}

\begin{algorithm}[t]
	\renewcommand{\thealgorithm}{}
	\caption{\hspace{-0.6em} \textbf{1} Numerical algorithm for solving the OCP (\ref{OCP1_cost})-(\ref{OCP1_sys}).} 
	\label{alg:Framwork}
	\begin{algorithmic}[1]
		\Require 
		Initialize $x_{t_0}$, $z_0=[u_{t_0}^{T},...,u_{t_{N_p}}^{T}]^{T}$;
		\Ensure  
		Optimal controller $u_{t_0}$, $t_0=0,...,N$;
		\For{$k=0$ \textbf{to} $N$}
		\For{$i=0$ \textbf{to} $\#step$}
		\State Calculate $\nabla J=\nabla J(x_{t_0},z_i)$; 
		\State Calculate $\nabla^2J=\nabla^2J(x_{t_0},z_i)$;
		\State $g_0(z_i) = (R+\nabla^2J)^{-1}\nabla J$;
		\For{$j=1$ \textbf{to} $i$}
		\State $g_j(z_i) = (R+\nabla^2J)^{-1}[\nabla J+Rg_{j-1}(z_i)]$;
		\EndFor
		\State $z_{i+1} = z_i - g_i(z_i)$;
		\EndFor
		\State Extract $u_{t_0}$ from $z_{i+1}$;
		\State Update initial state $x_{t_0}$;
		\EndFor
	\end{algorithmic}
\end{algorithm}

\subsection{Implementation}
In long-horizon OCP, performance is often compromised by numerical errors and model inaccuracies. MPC mitigates these issues by solving a sequence of local optimal control problems online. Instead of optimizing over the entire time horizon in one step, MPC focuses on a finite prediction horizon $k=t_0,...,t_{N_p}$ at each iteration, updating the control inputs based on the current system state. 

\par In practical implementation, solving the OCP (\ref{OCP1_cost})-(\ref{OCP1_sys}) at each sampling time is necessary within the MPC framework. During the optimization process, firstly, compute the gradient of the cost function with respect to the control variables using Lemma 1. Secondly, compute the hessian matrix of the cost function with respect to the control variables as described in Theorem 1. Finally, iteratively update the optimization algorithm by applying formulas (\ref{z_i+1_2}) through (\ref{g_M_1_2}). Repeat this process until the optimization converges. The details are summarized in Algorithm 1.

\section{Simulation and experimental results}
This section presents the simulation and experimental results to demonstrate the effectiveness of the proposed optimal control algorithm. First, the simulation validation of the LQR example is provided, followed by the simulation and experimental results of AGV trajectory tracking control.

\subsection{Analytical solution verification: LQR}
Consider the following linear discrete-time system:
\begin{align}
	\nonumber x_{k+1} = ax_k + bu_k,
\end{align}
where $x_k\in \mathbb{R}$ is the state variable with the given initial value $x_0$, $u_k\in \mathbb{R}$ the control variable, and $a$ and $b$ are the given constants. 
The associated quadratic cost function is given by
\begin{align}
	\nonumber J = \sum_{k=0}^{N-1}[qx_k^2 + ru_k^2] + px_N^2,
\end{align}
where $q\geq 0$, $r>0$, and $p\geq 0$. The optimal control problem outlined above can be analytically solved by solving the discrete Riccati equation to obtain the optimal controllers. This solution serves as a benchmark to assess the accuracy of Algorithm 1 in determining the controllers. 

\par
\begin{figure}[tbp]
	\centering
	\subfigure[\hspace*{0.0em}]{
		\includegraphics[width=0.23\textwidth]{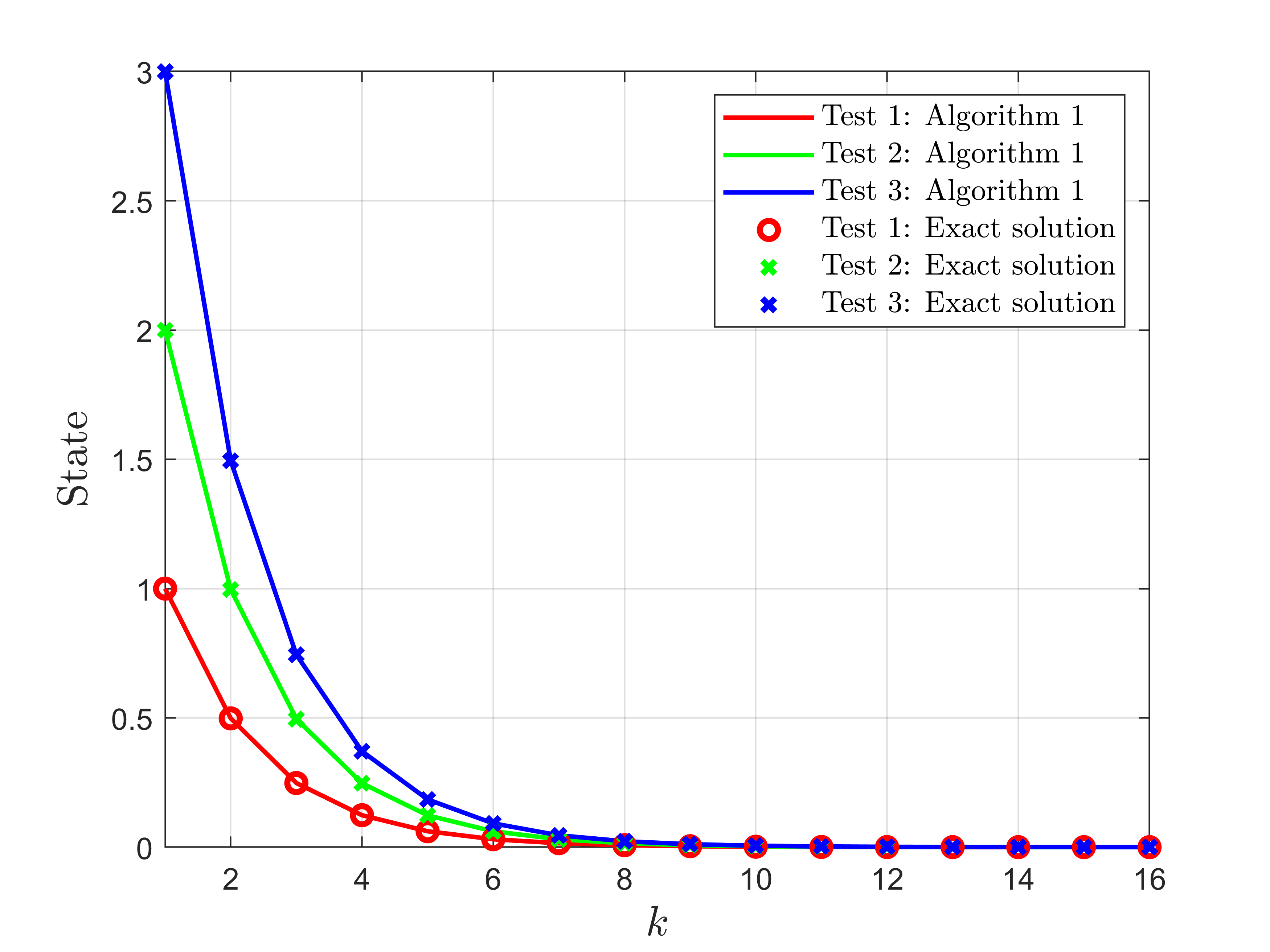}
	}
	\hspace*{-1em}
	\subfigure[\hspace*{0.0em}]{
		\includegraphics[width=0.23\textwidth]{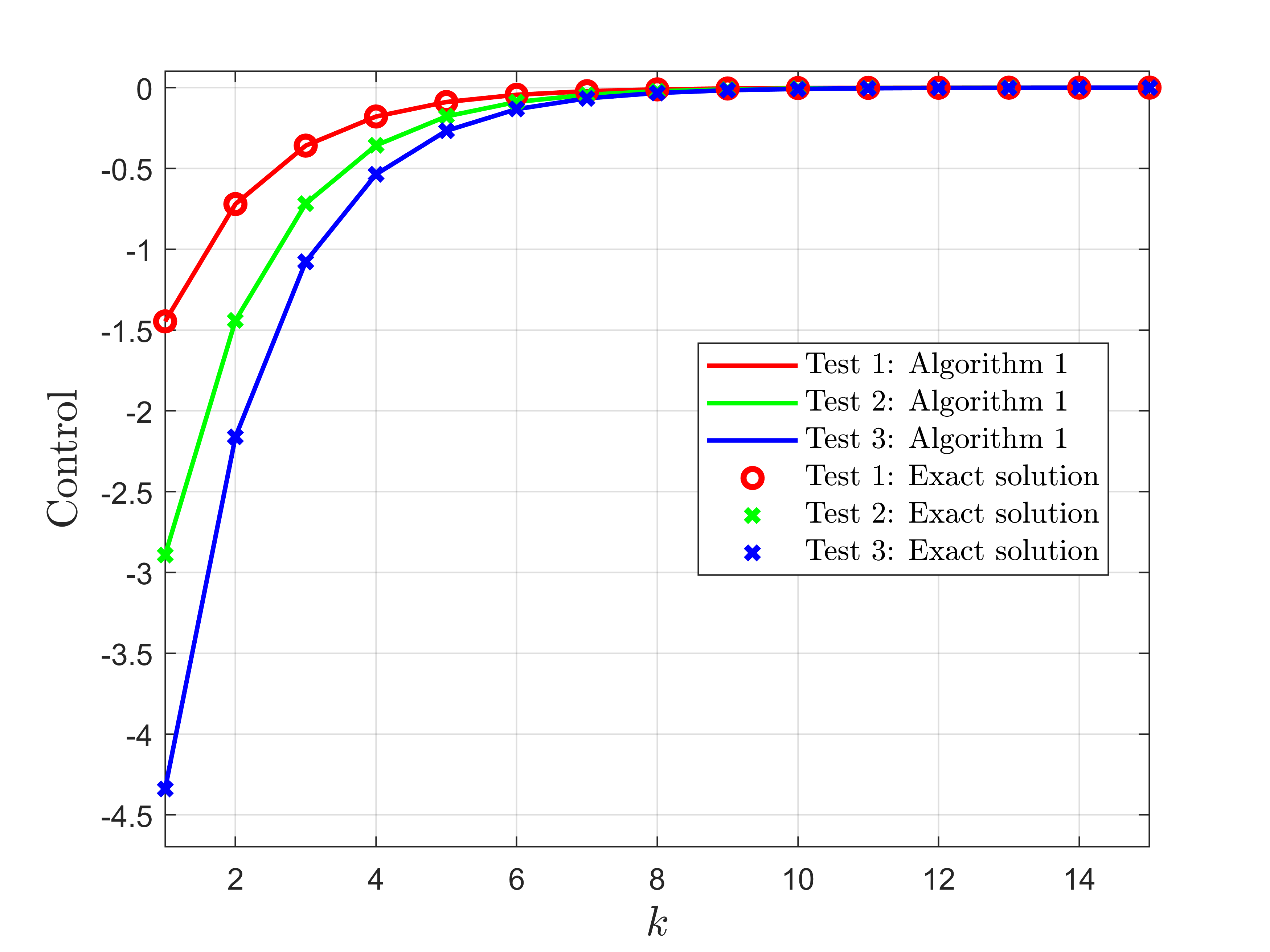}
	}
	\captionsetup{font={footnotesize}}
	\caption{Comparison of the optimal state and control obtained from Algorithm 1 and the analytical solution under three different initial conditions: $x_0=1.0, 2.0, 3.0$. (a) Optimal state. (b) Optimal control.} \label{fig_LQR}
\end{figure}

 In this example, we choose the parameters $a=1.8$, $b=0.9$, $q=1.0$, $r=3.0$, $p=3.0$, $N=15$, and $\bar R=0.1$. Besides, three independent simulations are conducted under different initial state value $x_0$, with Algorithm 1 running a single iteration of the outer loop.
Fig. \ref{fig_LQR} illustrates the convergence of optimal state and control obtained by Algorithm 1 and the analytical solution. It is evident that the results of Algorithm 1 nearly coincide with the analytical solution, as validated through three simulations. 

\subsection{Application: Trajectory tracking control of AGV}
\subsubsection{Problem Formulation}
Consider the following kinematic model of AGV
\begin{equation}
	\begin{cases}
		\dot x_t = v_t\cdot{\rm cos}\theta_t, \\
		\dot y_t = v_t\cdot{\rm sin}\theta_t,\\
		\dot \theta_t = \omega_t,
	\end{cases}                         \label{ODE_AGV}
\end{equation}
where $t\in[0,t_f]$.
The AGV's planar motion is characterized by three state variables: $x$, $y$, and $\theta$. Here, $x$ and $y$ denote the position in the plane, while $\theta$ represents the AGV's orientation. The control inputs are $v$ and $\omega$, which correspond to the linear and angular speeds of the AGV, respectively. 
\begin{figure*}[tbp]
	\centering
	\subfigure[\hspace*{0.0em}]{
		\includegraphics[width=0.33\textwidth]{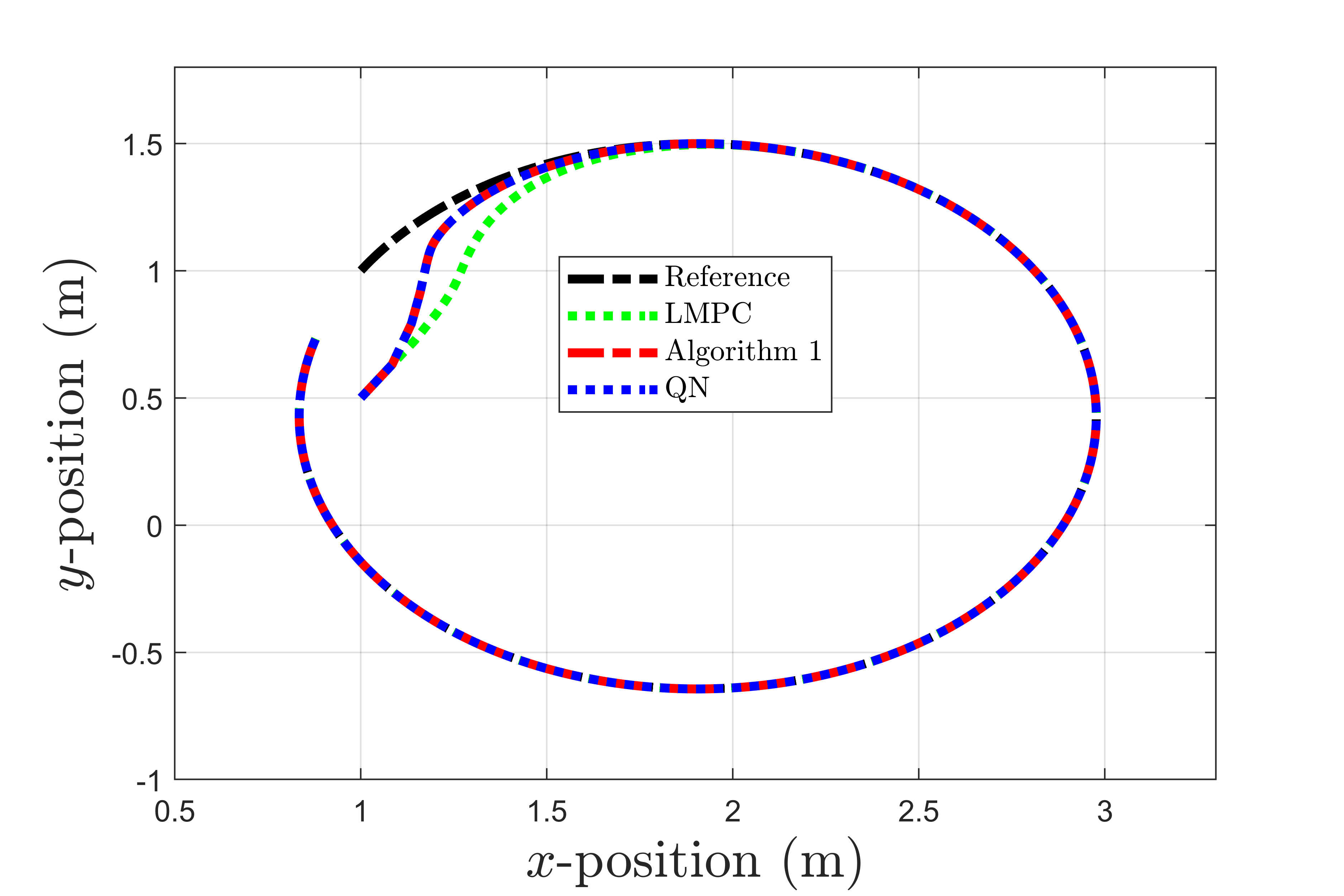}
	}
	\hspace*{-2.5em}
	\subfigure[\hspace*{0.0em}]{
		\includegraphics[width=0.33\textwidth]{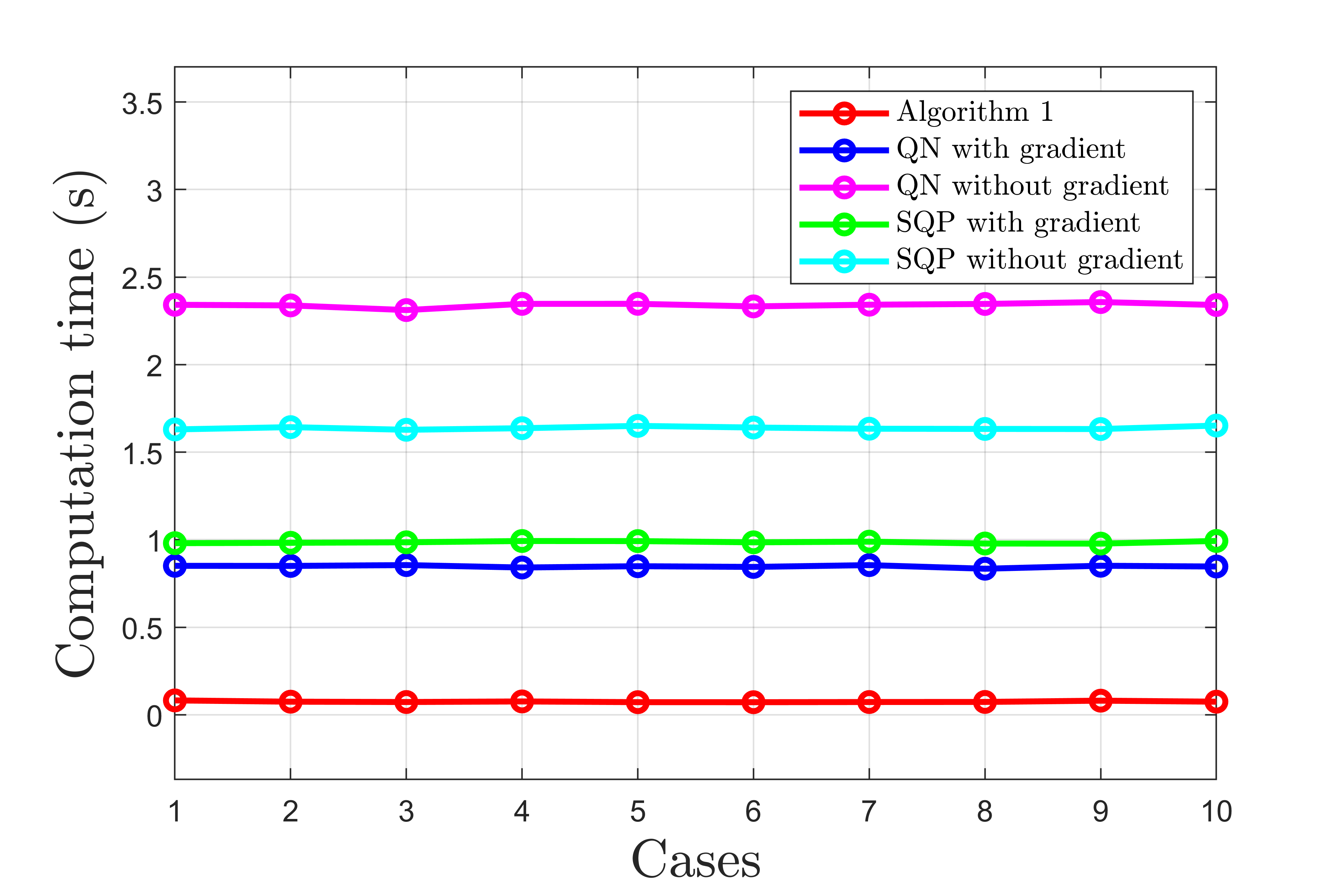}
	}
	\hspace*{-2.5em}
	\subfigure[\hspace*{0.0em}]{
		\includegraphics[width=0.33\textwidth]{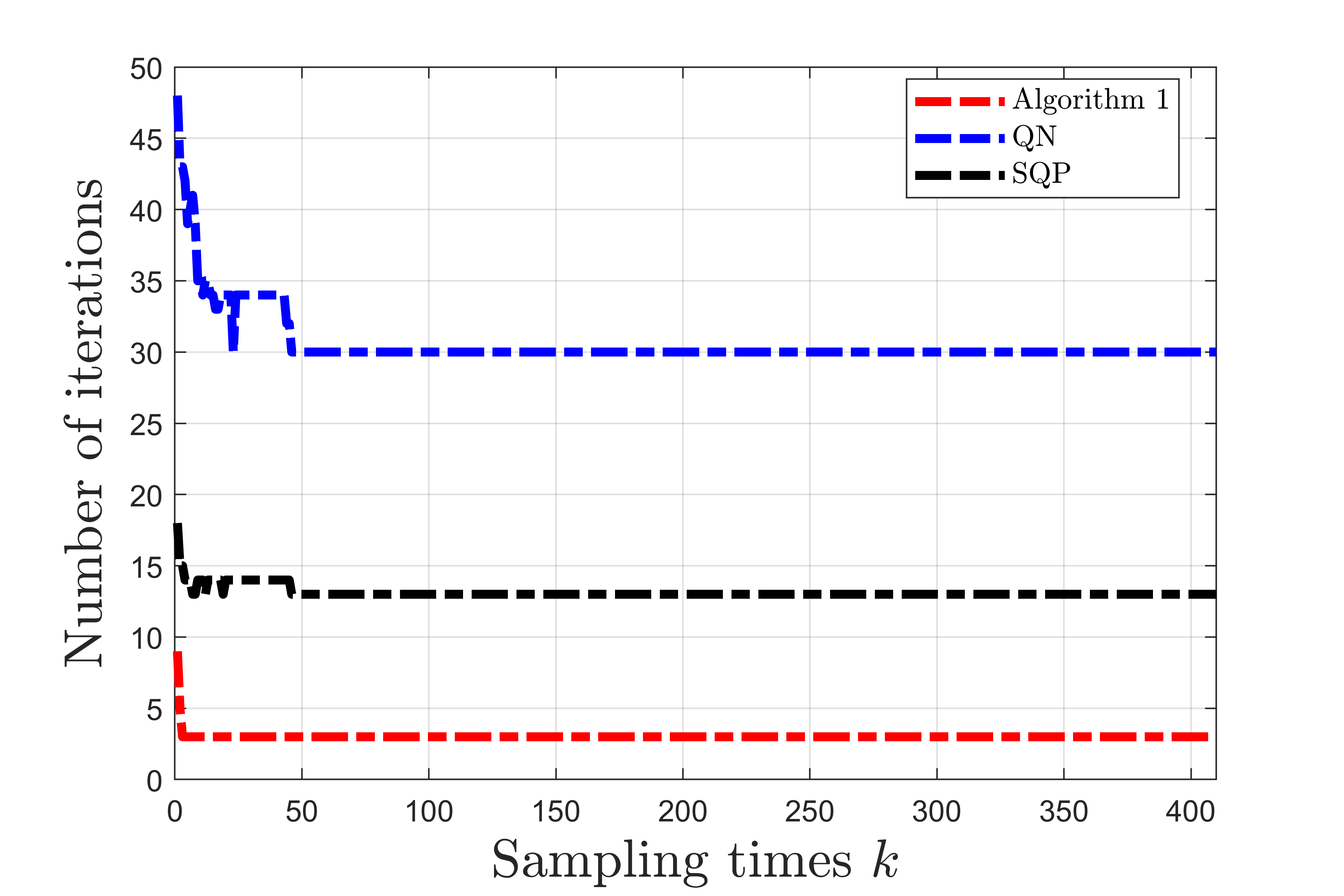}
	}
	\captionsetup{font={footnotesize}}
	\caption{Comparison of trajectory tracking results and relevant parameters among Algorithm 1, QN, and SQP. (a) The trajectory tracking results of Algorithm 1, LMPC and QN. (b) The computation time for running 10 independent simulations with each of the three algorithms. (c) The number of iterations required for optimization at each sampling time $k$ for the three algorithms.} \label{fig_AGV}
\end{figure*}
The system described by the equation (\ref{ODE_AGV}) is discretized using the forward Euler method, resulting in the following difference equations:
\begin{equation}
	\begin{cases}
		x_{k+1} = x_k + \Delta\cdot v_k{\rm cos}\theta_k , \\
		y_{k+1} = y_k + \Delta\cdot v_k{\rm sin}\theta_k,\\
		\theta_{k+1} = \theta_k + \Delta\cdot \omega_k,
	\end{cases}                         \label{AGV}
\end{equation}
where $k=0,1,...,N-1$ represents the sampling times, and $\Delta$ is the time step.
For convenience, let the state be denoted by $X_k=[x_k,y_k,\theta_k]^{T}$ and the control input by $U_k=[v_k,\omega_k]^{T}$. Additionally, denote the reference state and control input by $X_r=[x_r,y_r,\theta_r]^{T}$ and $U_r=[v_r,\omega_r]^{T}$, respectively. The error between the actual and reference values is given by $X_e=X_k-X_r$ and $U_e=U_k-U_r$.
The associated cost function is defined as
\begin{align}
	J = \sum_{k=0}^{N-1}(X_e^T\mathcal{Q}X_e + U_e^T\mathcal{R}U_e),  \label{AGV_cost}
\end{align} 
where $\mathcal{Q}$ and $\mathcal{R}$ are the weighting matrix. The objective is to design an optimal controller $U_k$ to enable the AGV to track the reference trajectory.

\subsubsection{Simulation Results}
\par In this example, set the terminal sampling time $N=410$, the time step $\Delta=0.05$, the prediction horizon $N_p=10$, the AGV's initial pose $X_0=[1.0,0.5,1.0]$, and the weight matrix $\mathcal{Q}={\text{diag}}[150,150,3]$, $\mathcal{R}={\text{diag}}[0.5,0.5]$, and $R=0.1I_{2N_p\times 2N_p}$. The simulations are conducted on a computer equipped with a 13th Gen Intel(R) Core(TM) i7-13620H processor and 16.0 GB of RAM, running Matlab R2020b. 

\par The OCP (\ref{AGV})-(\ref{AGV_cost}) can be solved using linearization MPC (LMPC) and existing optimization algorithms such as gradient descent, quasi-newton (QN), and sequential quadratic programming (SQP) \cite{cimen2004nonlinear, ren2020optimal, boyd2004convex}. Here, QN and SQP are selected as benchmarks for Algorithm 1. Although SQP is typically used for constrained optimization, it is also efficient for unconstrained problems. For implementation, QN and SQP algorithms are executed using Matlab's built-in solvers \texttt{fminunc} and \texttt{fmincon}, respectively. It should be noted that both QN and SQP can be implemented with or without gradient information.
Additionally, for each optimization, the initial values of the parameters are set to a zero vector, and the stopping criterion is defined as the infinity norm of $\nabla J$ being less than $10^{-6}$. 

\par Fig. \ref{fig_AGV} presents the tracking results for the circular trajectory along with a comparison of relevant metrics. In Fig. \ref{fig_AGV}(a), QN algorithm is used to validate the optimal controller obtained by Algorithm 1. The results demonstrate that the optimal state trajectory produced by Algorithm 1 aligns perfectly with the QN results, confirming the success of Algorithm 1 in solving the trajectory tracking problem. Besides, compared to LMPC, Algorithm 1 achieves faster response and higher precision control.
Fig. \ref{fig_AGV}(b) compares the computation time required by the three algorithms during the tracking process. When using QN and SQP algorithms, providing gradient information significantly reduces computation time, with both algorithms showing similar performance under these conditions. In contrast, Algorithm 1 achieves nearly a 10-fold reduction in computation time compared to these two algorithms. Fig. \ref{fig_AGV}(c) illustrates the number of iterations required to compute the optimal controller at each sampling time. It is clear that Algorithm 1 requires far fewer iterations than QN and SQP to reach the same level of accuracy.
\begin{figure}[t]
	\centering
	\includegraphics[width=0.47\textwidth]{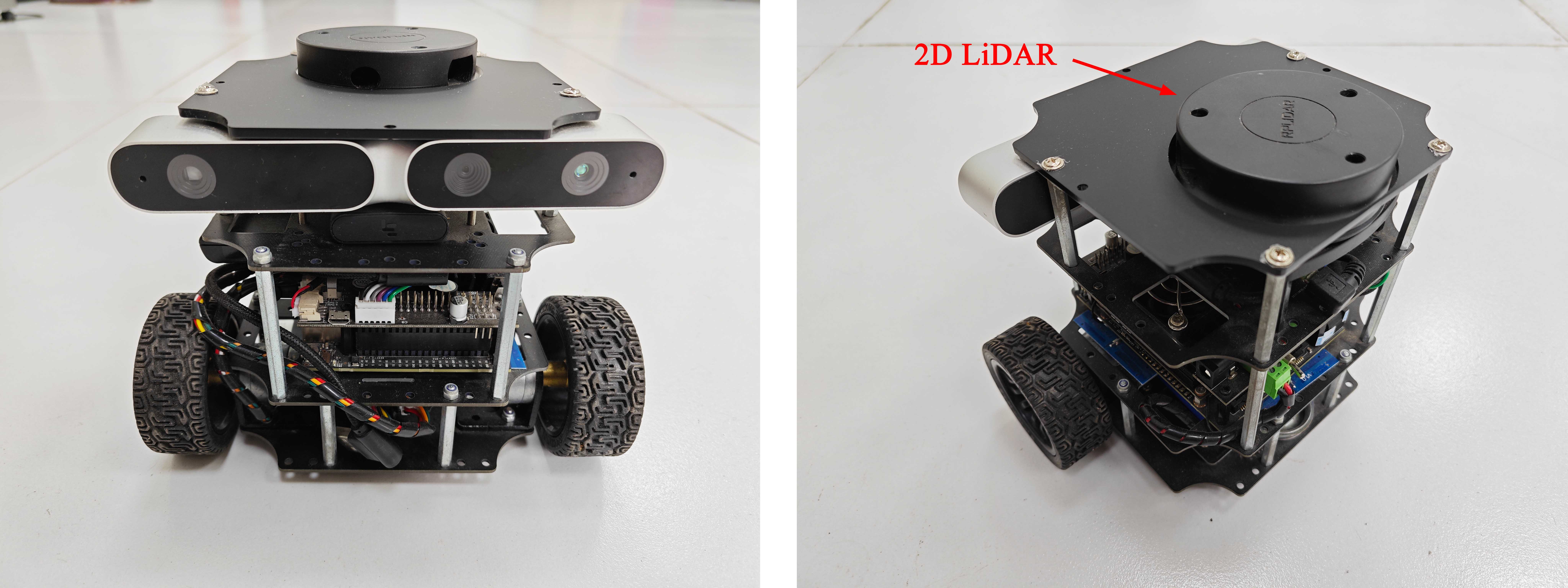}
	\captionsetup{font={footnotesize}}
	\caption{Experimental platform.} \label{car}
\end{figure}
\subsubsection{Experimental Results}
In the experiment, a closed indoor environment is selected. The AGV, as illustrated in Fig. \ref{car}, is primarily equipped with: 1) a 2D LiDAR for localization; 2) a Jetson Nano B01 mini PC to run Algorithm 1; and 3) a dual DC motor drive module AT8236 to control the motors. 

\par
\begin{figure}[tbp]
	\centering
	\subfigure[\hspace*{0.0em}]{
		\includegraphics[width=0.23\textwidth]{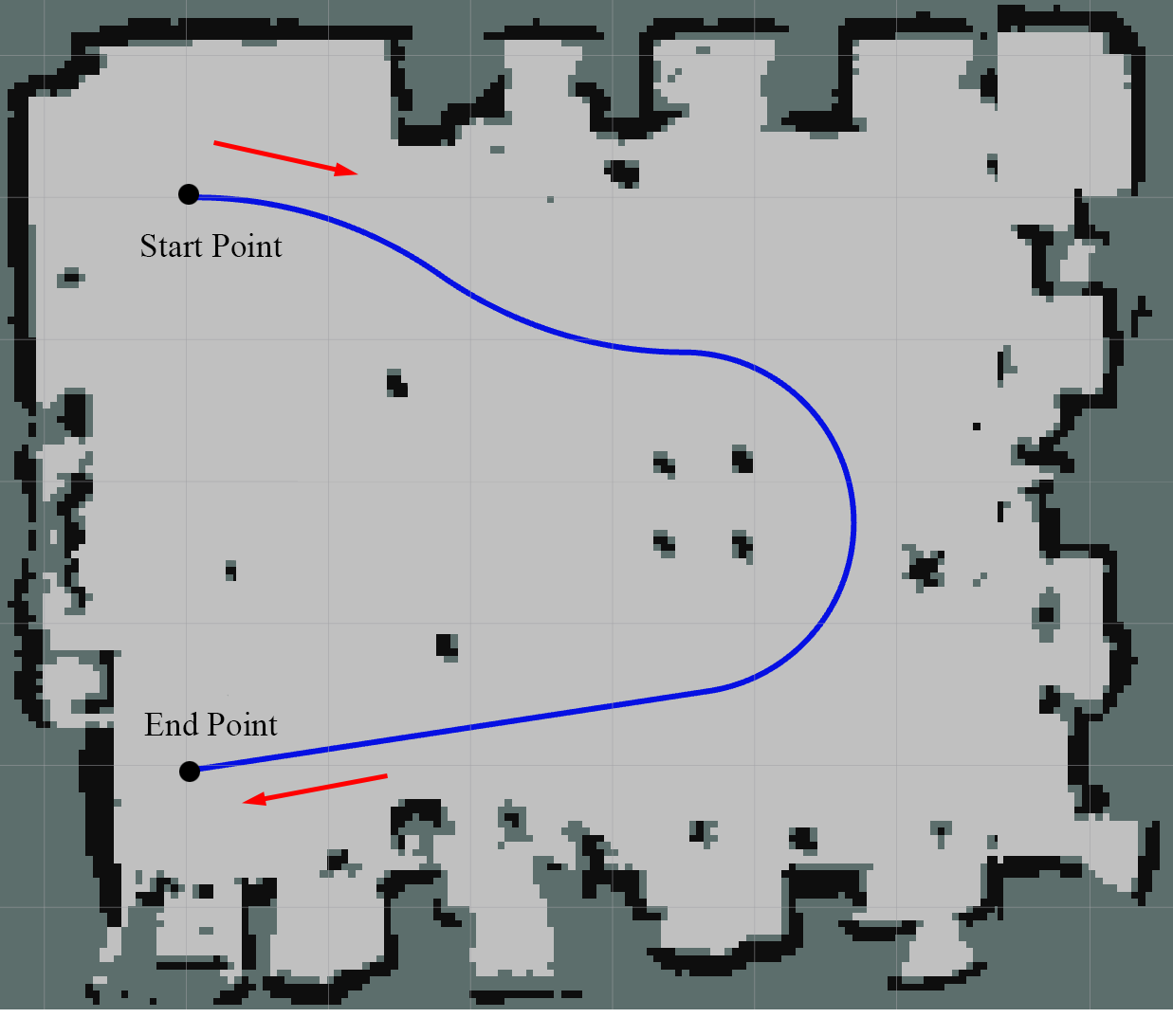}
	}
	\hspace*{-1em}
	\subfigure[\hspace*{0.0em}]{
		\includegraphics[width=0.23\textwidth]{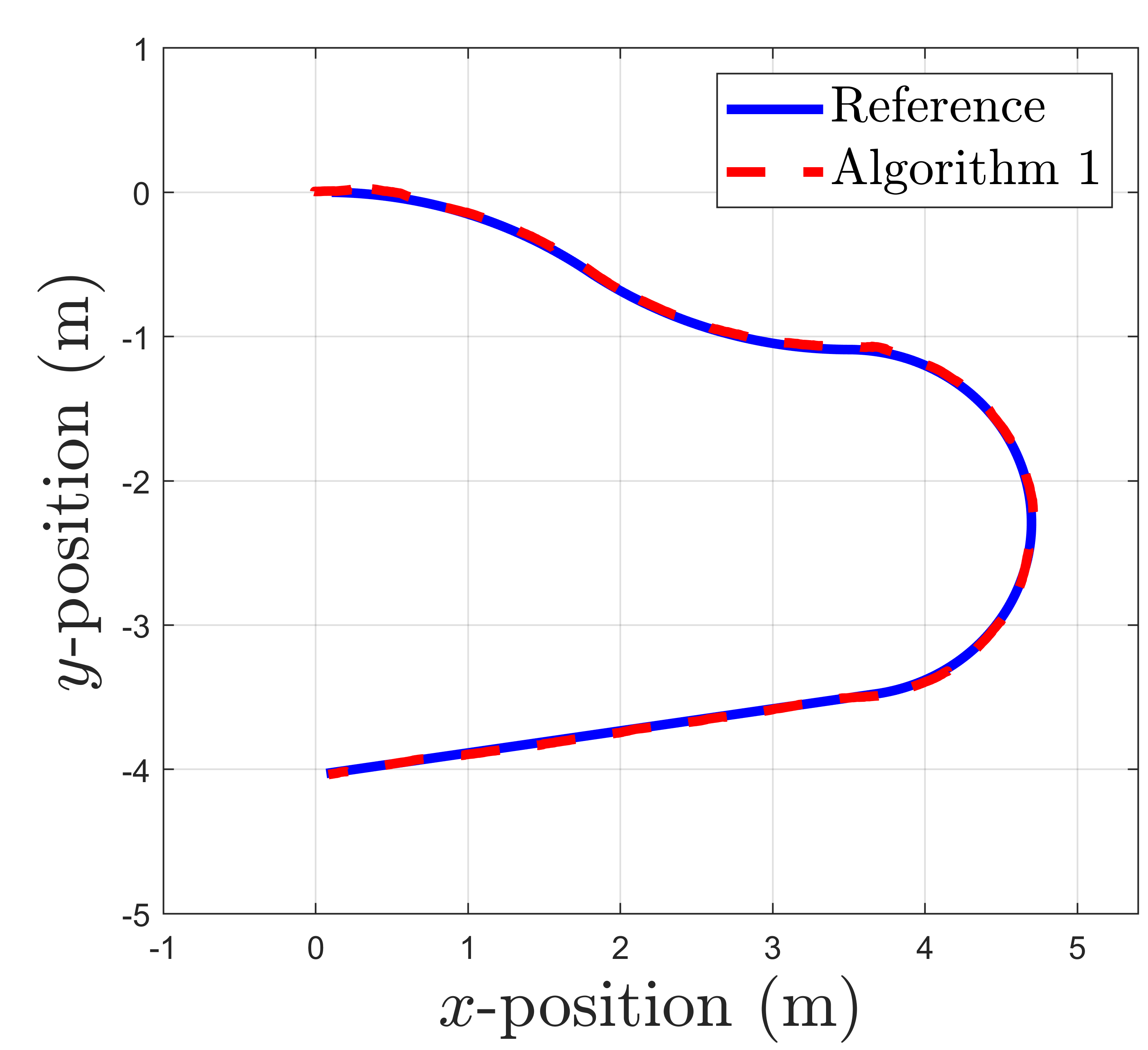}
	}
	\hspace*{-1em}
	\subfigure[\hspace*{0.0em}]{
		\includegraphics[width=0.53\textwidth]{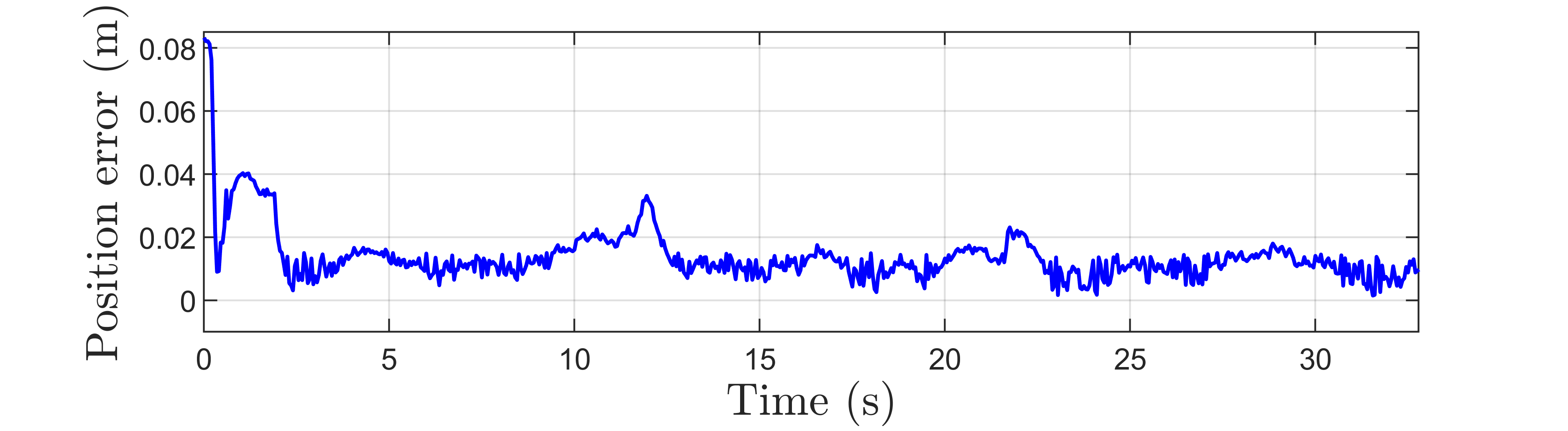}
	}
	\hspace*{-1em}
	\subfigure[\hspace*{0.0em}]{
		\includegraphics[width=0.53\textwidth]{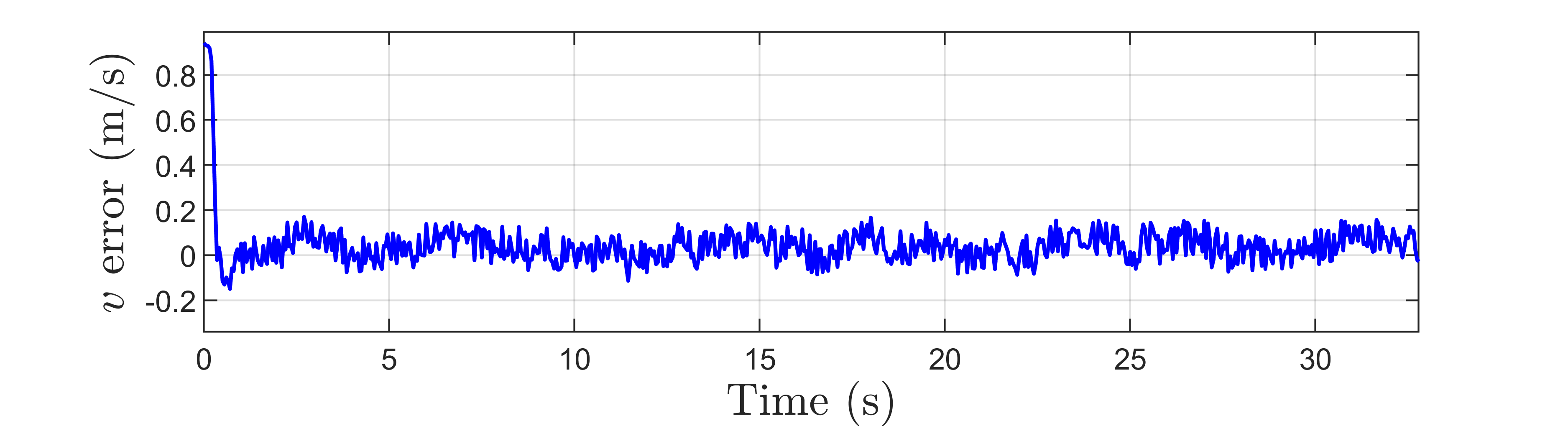}
	}
	\hspace*{-1em}
	\subfigure[\hspace*{0.0em}]{
		\includegraphics[width=0.53\textwidth]{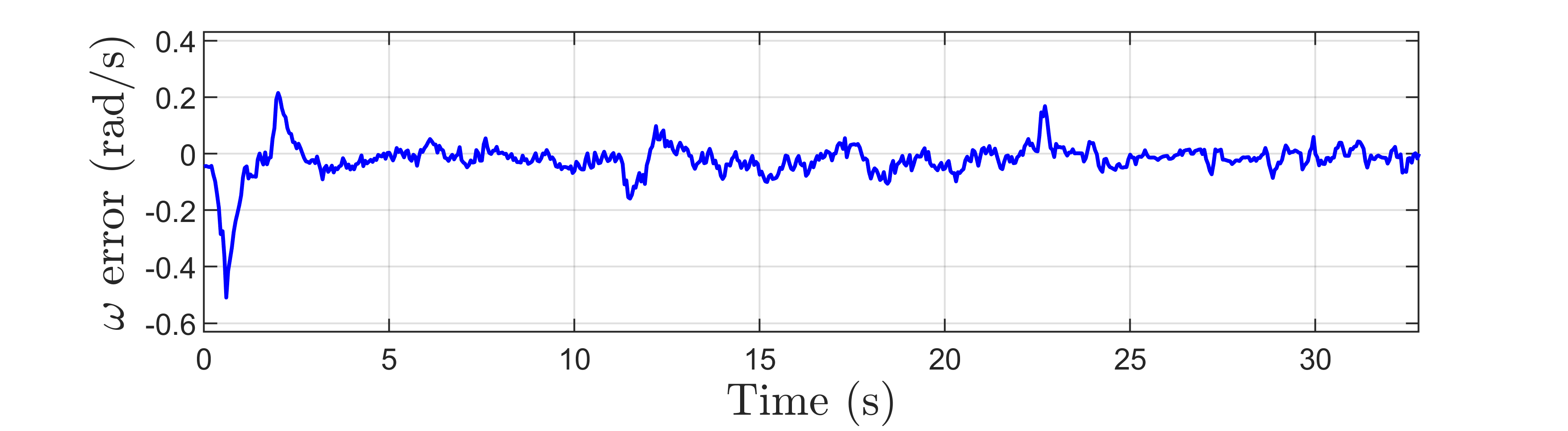}
	}
	\captionsetup{font={footnotesize}}
	\caption{The trajectory tracking results of AGV in indoor environment. (a) Map of the environment and the desired tracking trajectory. (b) AGV's tracking performance. (c) Trajectory tracking error. (d) Linear velocity $u$ error. (e) Angular velocity $\omega$ error.} \label{Experiment}
\end{figure}

 Fig. \ref{Experiment}(a) illustrates an indoor environmental map constructed using simultaneous localization and mapping (SLAM) technology, where the gray area denotes the accessible region for the robot. The blue line represents the tracking trajectory, with the starting point located at the top-left corner of the map and the endpoint at the bottom-left corner.
The validation process of the algorithm is carried out on the robot operating system (ROS) platform, with all the code written in C++. During the tracking procedure, we set up remote communication between a virtual machine running Ubuntu 20.04 and the vehicle to manage the AGV's starting, stopping, and data collection processes.

\par Fig. \ref{Experiment}(b) and \ref{Experiment}(c) illustrate the positional changes and error characteristics during the tracking process. The initial position of the vehicle is approximately 0.08 meters away from the starting point of the trajectory. At time points 12 seconds and 22 seconds, the vehicle exhibits slight fluctuations while turning. Apart from these three instances, the tracking error at other positions remains below 0.02 meters. Fig. \ref{Experiment}(d) and \ref{Experiment}(e) depict the variations in controller errors during the tracking process. Initially, the vehicle undergoes an acceleration phase to realign with the desired tracking trajectory, resulting in relatively large errors in linear and angular velocities during this phase. In the following stages, the vehicle keeps its speed close to the target values. There are slight increases in angular velocities during turns, while the variations at other positions remain relatively stable.

\section{Conclusion}
This paper presents a method that combines nonlinear optimization with the Pontryagin's maximum principle to solve discrete-time nonlinear OCP. The optimization objective is redefined from a new optimal control perspective, and a fast-converging iterative algorithm is designed, developing a novel approach for solving nonlinear OCP using optimal control theory. Additionally, explicit formulas based on FBDEs are derived to enhance the computational efficiency of evaluating gradient and second-order derivative for complex cost function. Finally, the proposed algorithm is applied to the trajectory tracking control of AGV, and its effectiveness is validated through simulations and experiments.


\bibliographystyle{Bibliography/IEEEtranTIE}
\bibliography{reference_with_linearization}

\end{document}